\documentclass{amsart}
\usepackage{amsthm,amsmath}
\usepackage{amscd,amssymb}
\usepackage{multicol}
\usepackage[pagebackref]{hyperref}

\newcommand{\0}{\mathcal O}
\newcommand{\cC}{\mathcal C}
\newcommand{\cD}{\mathcal D}
\newcommand{\cH}{\mathcal H}

\newcommand{\cS}{\mathcal S}
\newcommand{\cT}{\mathcal T}
\newcommand{\ba}{\widehat{a}}
\newcommand{\bb}{\widehat{b}}
\newcommand{\bE}{\widehat{E}}
\newcommand{\balpha}{\widehat{\alpha}}
\newcommand{\bh}{\widehat{h}}
\newcommand{\N}{\mathbb N}
\newcommand{\F}{\mathbb F}
\newcommand{\Z}{\mathbb Z}
\newcommand{\Q}{\mathbb Q}
\newcommand{\bQ}{\overline{\Q}}
\newcommand{\Qt}{{\mathbb Q}^\times}
\newcommand{\Qts}{\Q^{\times\,2}}
\newcommand{\eps}{\varepsilon}
\newfont{\cyr}{wncyb10}
\newcommand{\TS}{\mbox{\cyr Sh}}

\newcommand{\lra}{\longrightarrow}

\newcommand{\tors}{{\operatorname{tors}}}

\newcommand{\im}{{\operatorname{im}}\,}
\newcommand{\End}{{\operatorname{End}}}
\newcommand{\Sel}{{\operatorname{Sel}}}
\newcommand{\Cl}{{\operatorname{Cl}}}
\newcommand{\GL}{{\operatorname{GL}}}
\newcommand{\ts}{\textstyle}

\renewcommand{\Re}{\operatorname{Re}}

\newcounter{lemmacount}

\newtheorem{prop}[lemmacount]{Proposition}

\begin{document}

\author{Franz Lemmermeyer}
\title{Conics - a Poor Man's Elliptic Curves}
\maketitle
\tableofcontents
\vfill \eject

\section*{Introduction}
The aim of this article is to show that the arithmetic
of Pell conics admits a description which is completely analogous
to that of elliptic curves: there is a theory of $2$-descent
with associated Selmer and Tate-Shafarevich groups, and there
should be an analog of the conjecture of Birch and Swinnerton-Dyer.
For the history and a theory of the first $2$-descent, see
\cite{pell1,pell2,pell3}.
The idea that unit groups of number fields and the group of
rational points on elliptic curves are analogous is not new;
see e.g. \cite{D,DL,Lem,Zag} for some popularizations of this
point of view. It is our goal here to show that, for the case of 
the unit group of real quadratic number fields, this analogy can 
be made much more precise.

\section{The Group Law on Pell Conics and Elliptic Curves}

Let $F \in \Z[X,Y]$ a polynomial. If $\deg F = 2$, the affine curve
of genus $0$ defined by $F=0$ is called a conic. Let $d$ be a 
squarefree integer $\ne 1$ and define 
$$ \Delta = \begin{cases}
        d & \text{if}\ d \equiv 1 \bmod 4, \\
       4d & \text{if}\ d \equiv 2, 3 \bmod 4. \end{cases}$$
Then the curves $\cC: X^2 - \Delta Y^2 = 4$ are called Pell conics;
they are irreducible, nonsingular affine curves with a distinguished 
integral point $N = (2,0)$.

If $\deg F = 3$, the projective curve $E$ described by $F$ has genus
$1$ if it is nonsingular; if in addition it has a rational point,
then $E$ is called an elliptic curve defined over $\Q$. Elliptic
curves given by a Weierstra\ss{} equation $Y^2 = X^3 + aX + b$
are irreducible, nonsingular projective curves with a distinguished 
integral point $\0 = [0:1:0]$ at infinity.

Both types of curves have a long history: Pythagorean triples 
correspond to rational points on the Pell conic $X^4 + 4Y^2 = 4$, 
solutions of the Pell equations have been studied by the Greeks, 
the Indians, and the contemporaries of Fermat, such as Brouncker 
and Wallis. Problems leading to elliptic curves occur in the books 
of Diophantus and were studied by Bachet, Fermat, de Jonqui\`eres, 
Euler, Cauchy, Lucas, and Sylvester before Poincar\'e laid down his 
program for studying diophantine equations given by curves according 
to their genus.

\subsection{Group Law on Conics}
The group law on non-degenerate conics $C$ defined over a field $F$
is quite simple: fix any rational point $N$ on $C$; for finding
the sum of two rational points $A, B \in C(F)$, draw the line through
$N$ parallel to $AB$, and denote its second point of intersection with
$C$ by $A + B$. In the special case of Pell conics, the resulting
formulas can be simplified to

\begin{prop}\label{PGL}
Consider the conic $\cC: Y^2 - \Delta X^2 = 4$, and put $N = (2,0)$. 
Then the group law on $\cC$ with neutral element $N$ is given by
$$  (r,s) + (t,u) =  \Big(\frac{rt+\Delta su}2, \frac{ru + st}2\Big). $$
\end{prop}

It is now easily checked that the map sending points $(r,s) \in \cC(\Z)$ 
to the unit $\frac{r+s\sqrt{d}}2$ with norm $1$ in the quadratic number 
field $K = \Q(\sqrt{\Delta}\,)$ is a group homomorphism. Observe that 
the associativity of the geometric group law is equivalent to a special 
case of Pascal's theorem, which in turn is a very special case of 
Bezout's Theorem.

\subsection{Group Law on Elliptic curves}
Given an elliptic curve $E: y^2 = x^3 + ax + b$ defined over an 
algebraically closed field $K$, we define an addition law on $E$
by demanding that $A+B+C = 0$ for points $A, B, C \in E(K)$ if 
and only if $A, B, C$ are collinear. Since vertical lines intersect
$E$ only in two affine points, we have to regard $E$ as a projective
curve; then vertical lines intersect $E$ in two affine points as
well as in the point at infinity. Associativity follows geometrically
from a special case of Bezout's Theorem.

\section{The Group Structure}

Let us now compare the known results about the group structure of
Pell conics over the most common rings and fields. Generally, we 
will study conics in the affine plane over integral domains, and
elliptic curves in the projective plane over fields.

\subsection{Finite Fields}\label{Sp}

Let $\cC: x^2 - \Delta y^2 = 4$ be a Pell conic defined over
a finite field $\F_q$ with $q = p^f$ elements, and assume that 
$\cC$ is smooth, i.e. that $p \nmid \Delta$. Then
$$ \cC(\F_q) \simeq \Z/m\Z , \quad \text{where}\ 
             m = q - \Big(\frac{\Delta}{p}\Big)^f. $$ 
If $\Delta$ is a square mod $p$ and $p$ is odd, this is immediately
clear since there is an affine isomorphism between $\cC$ and
the hyperbolas $X^2 - Y^2 = 1$ and $XY = 1$; in particular, one 
has $\cC(\F_q) \simeq \F_q^\times = \GL_1(\F_q)$ in this case.

On the elliptic curve side, we know that 
$$E(\F_q) \simeq \Z/n_1\Z \oplus \Z/n_2\Z, n_2 \mid n_1, $$
Moreover, we have $\# E(\F_p) = (p+1) - a_p$, where $|a_p| \le 2\sqrt{p}$
by Hasse's theorem.

\subsection{$p$-adic Numbers}
If $p$ is an odd prime not dividing $\Delta$, then
$$ \cC(\Z_p) \simeq \begin{cases}
    \Z/(p-1) \oplus \Z_p & \text{if}\ (\frac{\Delta}{p}) = +1, \\
    \Z/(p+1) \oplus \Z_p & \text{if}\ (\frac{\Delta}{p}) = -1, \\
    \Z/2     \oplus \Z_p & \text{if}\ p \mid \Delta \ne -3, \\
    \Z/6     \oplus \Z_p & \text{if}\ p = 3, \ \Delta = -3.
   \end{cases}$$
Reduction modulo $p^k$ then yields 
$$  \cC(\Z/p^k) \simeq \begin{cases}
       \Z/(p-1) \oplus \Z/p^{k-1} & \text{if}\ (\frac{\Delta}{p}) = +1, \\
       \Z/(p+1) \oplus \Z/p^{k-1} & \text{if}\ (\frac{\Delta}{p}) = -1, \\
       \Z/2     \oplus \Z/p^{k}   & \text{if}\ p \mid \Delta \ne -3, \\
       \Z/6     \oplus \Z/3^{k-1} & \text{if}\ p = 3, \Delta = -3.
    \end{cases} $$

For elliptic curves $E/\Q_p$ we have a reduction map sending 
$\Q_p$-rational points to points defined over $\F_p$. 
The group $E_{ns}(\F_p)$ is the set of all nonsingular points 
of $E$ over $\F_p$. The subgroups $E_i(\Q_p)$ ($i = 0, 1$)
of $E(\Q_p)$ are defined as the inverse images of $E_{ns}(\F_p)$
and of the point of infinity of $E(\F_p)$ under the reduction map.
These groups sit inside the exact sequence
$$ \begin{CD} 
   0 @>>> E_1(\Q_p) @>>> E_0(\Q_p) @>>> E_{ns}(\F_p) @>>> 0.
   \end{CD} $$
The structure of $E_{ns}(\F_p)$ is known: if $E/\F_p$ is nonsingular,
it was discussed in Subsection \ref{Sp}; if $E/\F_p$ is singular, then
$E_{ns}(\F_p)$ is isomorphic to $\cC(\F_p)$ for a certain conic
$\cC$, and we say that $E$ has additive, multiplicative or split
multiplicative reduction if the conic is a parabola 
$(\cC(\F_p) \simeq \F_p$), a hyperbola ($\cC(\F_p) \simeq \F_p^\times$),
or an ellipse $(\cC(\F_p) \simeq \F_{p^2}[1]$, the group
of elements with norm $1$ in $\F_{p^2}$).

We also know hat $E_1(\Q_p) \simeq \Z_p$ and that the quotient group 
$E(\Q_p)/E_0(\Q_p)$ is finite. Its order $c_p$ is called the 
Tamagawa number for the prime $p$, and clearly $c_p = 1$ we have 
for all primes $p \nmid \Delta$. More exactly it can be shown
(albeit with some difficulty) that $c_p \le 4$ if $E$ has additive 
reduction, and that $c_p$ is the exact power of $p$ dividing $\Delta$ 
otherwise.

\subsection{Integral and Rational Points}
Now let us compare the structure of the groups of rational points:
for elliptic curves, we have the famous theorem of Mordell-Weil
that $E(\Q) \simeq E(\Q)_\tors \oplus \Z^r$, where 
$E(\Q)_\tors$ is the finite group of points of finite order,
and $r$ is the Mordell-Weil rank.  For conics, on the other hand, 
we have two possibilities: either $C(\Q) = \varnothing$
(for example if $C: x^2 + y^2 = 3$) or $C(\Q)$ is infinite, and
in fact not finitely generated (see Tan \cite{Tan}). The analogy 
can be saved, however, by looking at integers instead of rational
numbers: if $K$ is a number field with ring of $S$-integers $\0_S$, 
then
\begin{multicols}{2}
\begin{center}
$C(\0_S) \simeq C(\0_S)_\tors \oplus \Z^r $

$E(K) \simeq E(K)_\tors \oplus \Z^r$
\end{center}
\end{multicols}
\noindent
where $r \ge 0$ is called the Mordell-Weil rank. Shastri \cite{Shastri} 
computed the rank $r$ for the unit circle over number fields $K$
and $S = \varnothing$.

\medskip

Note that the group of integral points on the hyperbola $XY = 1$ 
is isomorphic to $R^\times = \GL_1(R)$. Number theoretic algorithms
working with the multiplicative group of $R = \Z/p\Z$ in general
have an analog for conics, as we will see in the next section.

\section{Applications}

\subsection{Primality Tests}

The classical primality test due to Lucas is the following:
\begin{prop}
An odd integer $n$ is prime if and only if there exists an
integer $a$ satisfying the following two conditions:
\begin{enumerate}
\item[i)] $a^{n-1} \equiv 1 \bmod n$;
\item[ii)] $a^{(n-1)/r} \not\equiv 1 \bmod n$ for every prime
           $r \mid (n-1)$.
\end{enumerate}
\end{prop}

This primality test is based on the multiplicative group 
$(\Z/n\Z)^\times$, that is, on the group $\cH(\Z/n\Z)$ of 
$\Z/n\Z$-rational points on the hyperbola $\cH: XY = 1$.
Something similar works for any Pell conic:

\begin{prop}\label{PrT}
Let $n \ge 5$ be an odd integer and $\cC: X^2 - \Delta Y^2 = 4$ 
a nondegenerate Pell conic defined over $\Z/n\Z$ with neutral
element $N = (2,0)$, and assume that $(\Delta/n) = -1$. Then $n$ 
is a prime if and only if there exists a point $P \in \cC(\Z/n\Z)$ 
such that
\begin{enumerate}
\item[i)] $(n+1)P = N$;
\item[ii)] $\frac{n+1}r\,P \ne N$ for any prime $r$ dividing $n+1$.
\end{enumerate}
\end{prop}

Of course, for both tests there are `Proth-versions' in which
only a part of $N \pm 1$ needs to be factored.

The following special case of Proposition \ref{PrT} is well known: 
if $n = 2^p-1$ is a Mersenne number, then $n \equiv 7 \bmod 12$ 
for $p \ge 3$, hence $(3/n) = -1$; if we choose the Pell conic
$\cC: X^2 - 12Y^2 = 4$ and and $P = (4,1)$, then the test above 
is nothing but the Lucas-Lehmer test. We remark in passing that 
Gross \cite{Gro} has come up with a primality test for Mersenne 
numbers based on elliptic curves.

\subsection{Factorization Methods}

The factorization method based on elliptic curves is very
well known. Can we replace the elliptic curve by conics?
Yes we can, and what we get is the $p-1$-factorization method
for integers $N$ if we consider the conic $\cH: xy = 1$, 
and some $p \pm 1$-factorization method for general Pell conics.
The details are easy to work out for anyone familiar with
Pollard's $p-1$-method.

\section{$2$-Descent}

Consider the Pell conic $\cC: X^2 - \Delta Y^2 = 4$.
Define a map $\alpha: \cC(\Q) \lra \Qt/\Qts$ by
$$ \alpha(x,y) = \begin{cases}
                (x+2)\Qts     & \text{if}\ x \ne -2, \\
                -\Delta \Qts  & \text{if}\ x = -2. \end{cases} $$
If $P = (x,y) \in \cC(\Z)$ with $x > 0$, then $P$ gives rise 
to an integral point on the descendant $\cT_a(\cC): aX^2 - bY2 = 4$, 
where $a$ is a positive squarefree integer determined by 
$\alpha(P) = a\Qts$, and $ab = \Delta$. Conversely, any integral 
point on some $\cT_a(\cC)$ gives rise to an integral point with 
positive $x$-coordinate on the Pell conic $\cC$.

It can be shown that $\alpha$ is a group homomorphism, and that
we have an exact sequence
$$ \begin{CD}
   0 @>>> 2\cC(\Z) @>>> \cC(\Z) @>{\alpha}>> \Qt/\Qts.
   \end{CD} $$
Moreover, we have $\# \im \alpha = 2^{r}$, where $r$ is the
Mordell-Weil-rank of $C(\Z)$, and the elements of $\im \alpha$
are represented by the first descendants $\cT_a$ with 
$\cT_a(\Z) \ne \varnothing$. Thus computing the Mordell-Weil 
rank is equivalent to counting the number of first descendants 
$\cT_a$ with an integral point (see \cite{pell3}).

The situation is completely analogous for elliptic curves
$E: Y^2 = X(X^2 + aX + b)$ with a rational point $(0,0)$ of 
order $2$, except that here we also have to consider the
$2$-isogenous curve $\bE: Y^2 = X(X^2 + \ba \,X + \bb)$,
where $\ba = -2a$ and $\bb = a^2 - 4b$. We have two Weil maps
$\alpha: E(\Q) \lra \Qt/\Qts$ and $\balpha: \bE(\Q) \lra \Qt/\Qts$,
and the Mordell-Weil rank is given by Tate's formula
$2^{r+2} = \# \im \alpha \cdot \# \im \balpha$. For more
information on the descent via $2$-isogenies we refer to 
Silverman \& Tate \cite{ST}.

\subsection{Selmer and Tate-Shafarevich Group}
The subset of descendants $\cT_a: ar^2 - bs^2 = 4$ 
with a rational point form a subgroup $\Sel_2(C)$ of $\Qt/\Qts$ 
called the $2$-Selmer group of $C$. Next we define the 
Tate-Shafarevich group $\TS_2(C)$ by the exact sequence
$$ \begin{CD} 
   1 @>>> \im \alpha @>>> \Sel_2(C) @>>> \TS_2(C) @>>> 1. 
   \end{CD} $$
In \cite{pell3} we have shown that the $2$-part of the  
Tate-Shafarevich group of the Pell conic $\cC: X^2 - \Delta Y^2 = 4$ 
is $\TS_2(\Z) \simeq \Cl^+(k)^2[2]$.

For a cohomological definition of Selmer and Tate-Shafarevich 
groups, we need to interpret conics as principal homogeneous 
spaces. Every conic $X^2 - \Delta Y2 = 4c$ is a principal 
homogeneous space for $\cC(\Q)$; this is to say that the map
$$\mu: \cD(\Z) \times \cC(\Z) \lra \cD(\Z):
\ts      \mu((u,v),(x,y)) = (\frac{ux + \Delta vy}2, \frac{vx+uy}2).$$
has the following properties:
\begin{enumerate}
\item $\mu(p,N) = p$ for all $p \in \cD(\bQ)$, where 
                      $N = (2,0)$ is the neutral element of $\cC$.
\item $\mu(\mu(p,P),Q) = \mu(p,P+Q)$ for all $p \in \cD(\bQ)$ and 
                      all $P, Q \in \cC(\bQ)$.
\item For all $p, q \in \cD(\Q)$ there is a unique $P \in \cC(\Q)$ 
      with $\mu(p,P) = q$.
\end{enumerate}
Here $\bQ$ denotes the algebraic closure of $\Q$. 

Note, however, that only those $\cD$ with $c \mid \Delta$ are 
principal homogeneous space for $\cC(\Z)$, i.e., satisfy the 
property that for all $p, q \in \cD(\Z)$ there is a $P \in \cC(\Z)$ 
with $\mu(p,P) = q$. Also observe that the conics $\cD$ with 
$c \mid \Delta$ can be written in the form $aX^2 - bY^2 = 4$
with $ab = \Delta$, that is, these are exactly the first 
descendants.

\subsection{Heights}
For a rational number $q = \frac{m}{n}$ in lowest terms, define
its height $H(q) = \log \max \{ |m|, |n|\}$; note that $H(0) = 0$ 
and $H(q) \ge 0$ for all $q \in \Q$. For rational points 
$P = (x,y) \in C(\Q)$ on a conic $C: X^2 - \Delta Y^2 = 4$ put
$H(P) = H(x)$.

Define the canonical height $\bh(P)$ by 
$$\bh(P) = \lim_{n \to \infty} \frac{H(2^nP)}{2^n}. $$
The canonical height $\bh$ on the Pell conic $\cC: X^2 - \Delta Y^2 = 4$
has all the suspected properties (and more):
\begin{enumerate}
\item $|\bh(P) - H(P)| < \log 4$;
\item $\bh(T) = 0$ if and only if $T \in \cC(\Q)_\tors$;
\item $\bh(mP) = m\bh(P)$ for all integers $m \ge 1$;
\item $\bh(P+Q) \le \bh(P) + \bh(Q)$;
\item the square of the canonical height satisfies the parallelogram
      equality 
      $$\bh(P-Q)^2 + \bh(P+Q)^2 = 2\bh(P)^2 + 2\bh(Q)^2$$
\end{enumerate}
for all $P, Q \in \cC(\Q)$.

In addition, there are explicit formulas for the canonical 
height. It is an easy exercise to show that every rational
point on a Pell conic has the form $P = (x,y)$ with 
$x = \frac{r}{n}$, $y = \frac{s}{n}$, and $(r,n) = (s,n) = 1$.
In this case we have
$$ \bh(P) = \begin{cases} 
          \log \frac{|r|+|s|\sqrt{\Delta}}2 & \text{if}\ \Delta > 0, \\
          \log |n|                          & \text{if}\ \Delta < 0.
     \end{cases} $$

The finiteness of $\cC(\Z_S)/2\cC(\Z_S)$ and the existence of a height
function implies the theorem of Mordell-Weil.

\section{Analytic Methods}

\subsection{Zeta Functions}

Both for conics and elliptic curves over $\Q$ there is an analytic
method that sometimes provides us with a generator for the group 
of integral or rational points on the curve. Before we can describe 
this method, we have to talk about zeta functions of curves.

Take a conic $C$ or an elliptic curve $E$ defined over the finite
field $\F_p$; let $N_r$ denote the cardinalities of the groups of
$\F_{p^r}$-rational points on $C$ and $E$ respectively, where we 
count solutions in the affine plane for $C$ and in the projective
plane for $E$. Then
$$ Z_p(T) = \exp\Big( \sum_{r=1}^\infty N_r \frac{T^r}r \Big) $$
is called the zeta function of $C$ or $E$ over $\F_p$.

For the parabola $C:y = x^2$, we clearly have $C(\F_q) \simeq \F_q$,
hence $N_r = p^r$, and we find
$$Z_p(T) = \exp \Big( \sum_{r=1}^\infty p^r \frac{T^r}r \Big)
     = \exp(- \log (1 -pT)) = \frac1{1-pT}.$$

For the conic $X^2 - \Delta Y^2 = 4$ we find after a little calculation 
$$ Z_p(T) = \frac1{(1-pT)(1-\chi(p)T)},$$
where $\chi$ is the Dirichlet character defined by $\chi(p) = (\Delta/p)$.
The substitution $T = p^{-s}$ turns this into
$$ \zeta_p(s; \cC) = \frac1{(1-p^{1-s})(1-\chi(p)p^{-s})}.$$

For nonsingular elliptic curves over $\F_p$ we similarly get
$$ Z_p(T) = \frac{P(T)}{(1-T)(1-pT)},$$
where $P(T) = qT^2 - a_pT + 1$ and $a_p$ is defined by 
$\# E(\F_p) = p+1 - a_p$.

\subsection{L-Functions for Conics}

Now we take the zeta function for each $p$ and multiply them 
together to get a global zeta function. The first factor
$1/(1-p^{1-s})$ gives us the product
$$ \prod_{p\, \text{odd prime}} \frac1{1-p^{1-s}} = \zeta(s-1)(1-2^{1-s}), $$
that is, essentially the Riemann zeta function.

The other factor, on the other hand, is more interesting:
$$ L(s,\chi) = \prod_p \frac1{1-\chi(p)p^{-s}} $$
is a Dirichlet $L$-function for the quadratic character 
$\chi = (\Delta/\,\cdot\,)$. This function converges on 
the right half plane $\Re s > 1$ and can be extended to 
a holomorphic function on the complex plane.

Now the nice thing discovered by Dirichlet (in his proof that
every arithmetic progression $ax+b$ with $(a,b)= 1$ contains
infinitely many primes) is that, for every nontrivial (quadratic) 
character $\chi$, $L(s,\chi)$ has a nonzero value at $s = 1$.
In fact, he was able to compute this value:

$$ L(1,\chi) = \begin{cases}
    h \cdot \frac{2\pi}{w\sqrt{|\Delta|}} & \text{if}\ \Delta < 0, \\
    h \cdot \frac{2\log \eps}{\sqrt{\Delta}} & \text{if}\ \Delta > 0
   \end{cases} $$
where $\chi(p) = (\Delta/p)$, and where $w$, $\Delta$, $h$ and $\eps> 1$
are the number of roots of unity, the discriminant, the class number
and the fundamental unit of $\Q(\sqrt{\Delta}\,)$.

The upshot is this: if $\Delta >0$, the group $\cC(\Z)$ has rank $1$;
by using only local information (numbers of $\F_{p^r}$-rational
points on $\cC$) we have constructed a function whose value
at $1$ gives, up to well understood constants, a power of
a generator of $\cC(\Z)$, namely the $h$-th power of the 
fundamental unit.

The functional equation of Dirichlet's $L$-function allows us to 
rewrite Dirichlet's formula as 
$$ \lim_{s \to 0} s^{-r} L(s,\chi) = \frac{2hR}{w}, $$
where $r = 0$ and $R = 1$ for $\Delta < 0$, and $r = 1$ 
and $R = \log \eps$ for $\Delta > 0$.

Observe that the evaluation of the $L$-funtion (which was defined
using purely local data) at $s = 0$ yields a generator of the
free part of the group $\cC(\Z)$ (which is a global object)!

\subsection{L-Functions for Elliptic Curves}
The really amazing thing is that exactly the same thing works
for elliptic curves of rank $1$: by counting the number $N_r$
of $\F_{p^r}$-rational points on $E$, we get a zeta function
$Z_p(T)$ that can be shown to have the form
$$ Z_p(T) = \frac{P(T)}{(1-T)(1-pT)}$$
for some polynomial $P(T) \in \Z[T]$ of degree $2$ (if $p$
does not divide the discriminant of $E$). In fact,
if $p \nmid E$ we have $P(T) = 1-a_pt + pt^2$, where   
$a_p = p + 1 - \# E(\F_p)$.

Put $L_p(s) = 1/P(p^{-s})$ and define the $L$-function
$$ L(s,E) = \prod_p L_p(s).$$ 
Hasse conjectured that this $L$-function can be extended 
analytically to the whole complex plane; moreover, there 
exists an $N \in \N$ such that
$$ \Lambda(s,E) = N^{s/2} (2\pi)^{-s} \Gamma(s) L(s,E) $$
satisfies the functional equation $\Lambda(s-2,E) = \pm \Lambda(s,E)$
for some choice of signs.
For curves with complex multiplication, this was proved by
Deuring; the general conjecture is a consequence of the now proved
Taniyama-Shimura conjecture.

\section{Birch--Swinnerton-Dyer}
\subsection{Birch and Swinnerton-Dyer for Elliptic Curves}
The conjecture of Birch and Swinnerton-Dyer for elliptic
curves predicts that $L(s,E)$ has a zero of order $r$ at $s =1$,
where $r$ is the rank of the Mordell-Weil group. More exactly, 
it is believed that 
$$ \lim_{s \to 1}\ \  (s-1)^r L(s;E) = 
    \frac{\Omega \cdot \# \TS(E/\Q) \cdot R(E/\Q) \cdot \prod c_p}
         {(\# E(\Q)_\tors)^2}, $$
where  $r$ is the Mordell-Weil rank of $E(\Q)$, $\Omega = c_\infty$ 
the real period, $\TS(E/\Q)$ the Tate-Shafarevich group, 
$R(E/\Q)$ the regulator of $E$ (some matrix whose entries are
canonical heights of basis elements of the free part of $E(\Q)$),
$c_p$ the Tamagawa number for the prime $p$ (trivial for all primes
not dividing the discriminant), and $E(\Q)_\tors$ the torsion group of $E$.

\subsection{Birch and Swinnerton-Dyer for Conics}

We now want to interpret Dirichlet's class number formula 
in a similar way. Let $k = \Q(\sqrt{\Delta}\,)$ denote the
quadratic number field associated to the Pell conic
$\cC: X^2 - \Delta Y^2 = 4$. Then we conjecture that 
there is a cohomological definition of the Tate-Shafarevich group
$\TS(\cC)$ whose $2$-torsion coincides with the group $\TS_2(\cC)$
defined above, and that we have 
$$ \TS(\cC) \simeq \Cl^+(k)^2. $$ 
If we (preliminarily) define the Tamagawa numbers by 
$$ c_p = \begin{cases}  
         2 & \text{if}\ p \mid \Delta, \\
         1 & \text{otherwise,} \end{cases} $$
then Gauss's genus theory implies that 
$$ \prod c_p = 2(\Cl^+(k):\Cl^+(k)^2).$$

Thus if we put $\Omega = \frac12$, then    
$\Omega \cdot \# \TS(\cC) \cdot \prod c_p = h^+$ 
equals the class number of $k$ in the strict sense, 
hence is equal to $2^u \cdot h$, where $u = 1$ if $N \eps = +1$, 
and $u = 0$ otherwise.

If $\Delta > 0$, let $\eta > 1$ denote a generator of the free 
part of $\cC(\Z)$; then the regulator of $\cC$ equals 
$\bh(\eta) = \log \eta$. Now we find $R(\cC) = 2^{1-u}R$, hence 
$\Omega \cdot \# \TS(\cC) \cdot R(\cC) \cdot \prod c_p 
  = h^+ \log \eta = 2 h R$; this also holds for $\Delta < 0$
if we put $R = 1$.
  
Finally, $\cC(\Z)_\tors$ is the group of roots of unity 
contained in $k$, and we find
$$ \frac{2hR}{w} = 
  \frac{\Omega \cdot \# \TS(\cC) \cdot R(\cC) \cdot \prod c_p}
       {\# \cC(\Z)_\tors} $$ 
in (almost) perfect analogy to the Birch--Swinnerton-Dyer 
conjecture for elliptic curves.

In fact, the analogy would be even closer if we would
replace $\# \cC(\Z)_\tors$ by $(\# \cC(\Z)_\tors)^2$ and
adjust the formulas for $c_2$ and $c_3$ for the two
Pell conics with nontrivial torsion; this would also 
allow us to put $\Omega = 1$. 

\section{Summary}

The analogy between Pell conics and elliptic curves is
summarized in the following table:

\bigskip
\begin{tabular}{l||l|l|l} 
 		   & $\GL_1$  & Pell conics & elliptic curves \\ 
                      \hline \hline
group structure on & affine line & affine plane  & projective plane \\
defined over       & rings    & rings         & fields \\
group elements     & $S$-units & $S$-integral points & rational points \\
group structure    & $\Z/2 \oplus \Z^{\#S}$  
                   & $C(\Z_S)_\tors \oplus \Z^r$ 
                   & $E(\Q)_\tors \oplus \Z^r$ \\
associativity      & clear & Pascal's Theorem & Bezout's Theorem \\
factorization alg. & $p-1$ & $p \pm 1$ & ECM \\
primality tests    & Lucas-Proth  & Lucas-Lehmer & ECPP \\
\TS                & 1  & $\Cl^+(k)^2$  &  ?  \\
L-series           & $\Z$ & quadratic field & modular form            
\end{tabular}
\bigskip

Moreover, cyclotomic fields are for Pell conics what modular curves 
are for elliptic curves, and cyclotomic units correspond to Heegner 
points. The analog of Heegner's Lemma (if a curve of genus $1$ of
the form $Y^2 = f_4(X)$, where $f_4$ is a quartic polynomial with 
rational coefficients, has a $K$-rational point for some number field 
$K$ of odd degree, then the curve has a rational point; cf. \cite{Hee}) 
is due to Nagell \cite{Nag}, who proved the same result with $f_4$ 
replaced by a quadratic polynomial $f_2$.

\section{Questions}
Although the arithmetic of conics is generally regarded as
being almost trivial, there are a lot of questions that are
still open. The main problem is a good definition of the
Tamagawa numbers in the case of conics, a cohomological
description of the Selmer and Tate-Shafarevich groups, and
the proof of $\TS(\cC) \simeq \Cl^+(k)^2$. 

The next problem is the analytic construction of generators
of $\cC(\Z_S)$ if $S \ne \varnothing$. This suggests looking
at the  Stark conjectures, which predict that we can construct 
certain units (actually $S$-units) in number fields. It seems,
however, that we cannot hope to find ``independent'' elements
(see \cite{Tang}). 

On a simpler level there's the question whether iterated
$2$-descents on Pell conics provide an algorithm for 
computing the fundamental unit that is faster than current
methods. And how does $3$-descent on Pell conics work? 

We can also think of generalizing the approach described here:
the groups $\GL_1$ and the Pell conics are special norm tori 
in the theory of algebraic groups, and there's the question of
how much of the above carries over to the more general situation.
The norm-$1$ tori associated to pure cubic fields can be
described geometrically as cubic surfaces $\cS$; do the groups of 
integral points on $\cS$ admit a geometric group law?
It is known that the groups of rational points on cubic surfaces
coming from norm forms satisfy the Hasse principle; is there
a connection between the $3$-class groups of these fields and 
the Tate-Shafarevich groups on $\cS$ defined as above as the 
obstruction to lifting the Hasse principle from rational to 
integral points?

On the elliptic curve side, there are a few questions 
suggested by the analogy worked out in this article. For example,
is there a natural group whose order equals $\# \TS(E) \cdot \prod c_p$?
Recall that $\exp(\bh(P))$ is algebraic for rational points on
Pell conics; are there meromorphic functions $F$ such that 
$F(\bh(P))$ is algebraic for rational points $P$ on elliptic curves,
at least for curves with complex multiplication?

\section*{Acknowledgments}
This article owes a lot to work done while I was at
the University of Seoul in August 2002; I would like to 
thank Soun-Hi Kwon for the invitation and the hospitality.

\end{document}